\documentclass[10pt]{amsart}

\usepackage{amssymb}
\usepackage{enumitem}

\usepackage[english]{babel}

\newtheorem{theorem}{Theorem}[section]

\newtheorem{e-proposition}[theorem]{Proposition}
\newtheorem{corollary}[theorem]{Corollary}
\newtheorem{e-definition}[theorem]{Definition\rm}

\newtheorem{example}[theorem]{Example}

\newcommand{\bbQ}{{\mathbb Q}}
\newcommand{\bbR}{{\mathbb R}}

\newcommand{\bbZ}{{\mathbb Z}}

\newcommand{\bfH}{{\mathbf H}}

\newcommand{\calO}{\mathcal{O}}

\newcommand{\pr}{\operatorname{pr}}

\newcommand{\Ker}{\operatorname{Ker}}

\newcommand{\Homeo}{\operatorname{Homeo}}

\newcommand{\SL}{\operatorname{SL}}
\newcommand{\SO}{\operatorname{SO}}

\newcommand{\Aut}{\operatorname{Aut}}

\newcommand{\PSL}{\operatorname{PSL}}
\newcommand{\PGL}{\operatorname{PGL}}
\newcommand{\GL}{\operatorname{GL}}

\newcommand{\QI}{\operatorname{QI}}

\begin{document}

\title{On the structure and arithmeticity of lattice envelopes}

\author{Uri Bader}
\address{Technion, Haifa}
\email{uri.bader@gmail.com}

\author{Alex Furman}
\address{University of Illinois at Chicago, Chicago}
\email{furman@math.uic.edu}

\author{Roman Sauer}
\address{Karlsruhe Institute of Technology}
\email{roman.sauer@kit.edu}


\maketitle

\begin{abstract}
We announce results about the structure and arithmeticity of all possible lattice embeddings 
of a class of countable groups which
encompasses all linear groups with simple Zariski closure, all groups with non-vanishing first $\ell^2$-Betti number,
word hyperbolic groups, and, more general, convergence groups.
\end{abstract}

\section{Introduction}
\label{sec:intro}

Let $\Gamma$ be a countable group. We are concerned with the study of its 
lattice envelopes, i.e.~the locally compact groups containing~$\Gamma$ 
as a lattice. 
We aim at structural results that impose no restrictions on the ambient
locally compact group and only abstract group-theoretic conditions on~$\Gamma$. 
We say that $\Gamma$ satisfies~$(\dagger)$ if every finite index subgroup of 
a quotient of $\Gamma$ by a finite normal subgroup 
\begin{enumerate}[label=($\dagger$\arabic*)]
	\item\label{irred} is not isomorphic to a product of two infinite groups, and 
	\item\label{caf} does not possess infinite amenable commensurated subgroups, and
	\item\label{nbc} satisfies: For a normal subgroup $N$ and a commensurated
	subgroup $M$ with $N\cap M=\{1\}$ there exists a finite index subgroup $M'<M$
	such that $N$ and $M'$ commute.
\end{enumerate}

The relevance of~\ref{irred} should be clear, the relevance of~\ref{caf} is that it
yields an information about all possible lattice envelopes of $\Gamma$~\cite{envelopes}:

\begin{e-proposition}\label{prop: compact amenable radical}
	Let $\Gamma$ be a lattice in a locally compact group $G$. If $\Gamma$ has no infinite amenable commensurated
	subgroups, then the amenable radical $R(G)$ of $G$ is compact.
\end{e-proposition}

The role of~\ref{nbc} is less transparent, but be aware of the obvious observation: If $M,N$ are both normal
with $N\cap M=\{1\}$, then $M$ and $N$ commute. There are lattices $\Gamma$ in
$\SL_n(\bbR)\times \Aut(T)$, where $T$ is the universal cover of the $1$-skeleton $B^{(1)}$ of the Bruhat-Tits
building of $\SL_n(\bbQ_p)$ (see~\cite[6.C]{caprace+monod-discrete} 
and~\cite[Prop.~1.8]{burger+mozes-localglobal}). They are built in such a way that $N:=\pi_1(B^{(1)})$, which is a free
group of infinite rank, is a normal subgroup of $\Gamma$. Let $U<\Aut(T)$ be the stabilizer of a vertex.
Then $M:=\Gamma\cap (\SL_n(\bbR)\times U)$ is commensurated, but $M,N<\Gamma$ violate~\ref{nbc}. Moreover, this group
$\Gamma$ satisfies~\ref{irred} and~\ref{caf}. 
\begin{e-proposition} \label{examples}
Linear groups with semi-simple Zariski closure satisfy the conditions \ref{caf} and \ref{nbc}.
Groups with some positive $\ell^2$-Betti number satisfy the condition \ref{caf}.
All the $(\dagger)$ conditions are satisfied by
all linear groups with simple Zariski closure,
by all groups with positive first $\ell^2$-Betti number,
by all non-elementary word hyperbolic or, more generally, convergence groups.
\end{e-proposition}

For a concise formulation of our main result, we introduce the following notion of $S$-arithmetic lattice embeddings
up to tree extension: Let $K$ be a number field. Let $\bfH$ be a connected, absolutely simple adjoint $K$-group, and let
$S$ be a set of (equivalence classes) of places of $K$ that contains every infinite place for which
$\bfH$ is isotropic and at least one finite place for which $\bfH$ is isotropic. Let $\calO_S\subset K$
denote the $S$-integers. The (diagonal) inclusion of (a finite index subgroup of) $\bfH(\calO_S)$
into $\prod_{\nu\in S}\bfH(K_\nu)^+$
is the prototype of an $S$-arithmetic lattice. 
Let $H$ be a group obtained from
$\prod_{\nu\in S}\bfH(K_\nu)^+$ by possibly replacing each factor $\bfH(K_\nu)$ with $K_\nu$-rank~$1$ by
an intermediate closed subgroup $\bfH(K_\nu)^+<D<\Aut(T)$ where $T$ is a tree with a cofinite $\bfH(K_\nu)^+$-action.
The lattice embedding $\bfH(\calO_S)$ into $H$ is called an \emph{$S$-arithmetic lattice embedding
up to tree extension}. \smallskip

A typical example is $\SL_2(\bbZ[1/p])$ embedded diagonally as a lattice into $\SL_2(\bbR)\times \SL_2(\bbQ_p)$. The
latter is a closed cocompact subgroup of $\SL_2(\bbR)\times\Aut(T_{p+1})$, where $T$ is the Bruhat-Tits tree of
$\SL_2(\bbQ_p)$, i.e.~a $(p+1)$-regular tree. So $\SL_2(\bbZ[1/p])<\SL_2(\bbR)\times\Aut(T_{p+1})$ is an
$S$-arithmetic lattice embedding up to tree extension. We now state the 
main result~\cite{envelopes}: 

\begin{theorem}\label{thm: main}
	Let $\Gamma$ be a finitely generated group satisfying $(\dagger)$, e.g one of the groups considered in Proposition~\ref{examples}. 
	Then every embedding
	of $\Gamma$ as a lattice into a locally compact group $G$ is, up to passage to finite index subgroups and dividing
	out a normal compact subgroup of $G$, isomorphic to one of the following cases:
	\begin{enumerate}
		\item an irreducible lattice in a center-free, semi-simple Lie group without compact factors;
		\item an $S$-arithmetic lattice embedding up to tree extension;
		\item a lattice in a totally disconnected group with trivial amenable radical.
	\end{enumerate}
	The same conclusion holds true if one replaces the assumption that $\Gamma$ is finitely generated by
	the assumption that $G$ is compactly generated.
\end{theorem}

Finite generation of $\Gamma$ implies compact generation of any locally compact group
containing $\Gamma$ as a lattice.
The examples above for $n\ge 3$ show that condition~\ref{nbc} in Theorem~\ref{thm: main} is indispensable. 
Since non-uniform lattices with a uniform upper bound on the order of finite subgroups do not exist in totally disconnected groups, our main theorem yields the following classification of non-uniform lattice embeddings.

\begin{corollary} 
	Let $\Gamma$ be a group that satisfies $(\dagger)$ and admits a uniform upper bound on the order of all finite subgroups. Then every non-uniform lattice embedding of $\Gamma$ into a compactly generated
	locally compact group $G$ is, up to passage to finite index subgroups and dividing
	out a normal compact subgroup of $G$, either a lattice in a center-free, semi-simple Lie group without compact factors
	or an $S$-arithmetic lattice embedding up to tree extension.
\end{corollary}

The following arithmeticity theorem~\cite{envelopes} is at the core of the proof of Theorem~\ref{thm: main}. Actually,
it is a more general version that is used in which we drop condition~\ref{irred} (see the comment
in~Step~3 of Section~\ref{sec: ideas}).

In the proof of Theorem~\ref{thm: main} we only need Theorem~\ref{thm: adelic arithmeticity}
in the case where $D$, thus $L\times D$, is compactly generated which means that the set $S$
of primes is finite. Caprace-Monod~\cite[Theorem~5.20]{caprace+monod-discrete}
show Theorem~\ref{thm: adelic arithmeticity}
for compactly generated~$D$ and under the hypothesis that $L$ is
the $k$-points of a simple $k$-group (where $k$ is a local field) 
but the latter hypothesis is too restrictive
for our purposes. Morevover, our proof of Theorem~\ref{thm: adelic arithmeticity} 
does not become much easier if we assume compact generation from the beginning. 
Regardless of its role in Theorem~\ref{thm: main} we consider the following result as a first step in the classification
of lattices in locally compact groups that are not necessarily compactly generated.

\begin{theorem}\label{thm: adelic arithmeticity}
		Let $L$ be a connected center-free semi-simple Lie group
		without compact factors, and let
		$D$ be a totally disconnected locally compact group without compact normal subgroups.  
		Let $\Gamma<L\times D$ be a lattice such that the projections of $\Gamma$ to both $L$ and $D$ are dense
		and the projection of $\Gamma$ to $L$ is injective and $\Gamma$ satisfies~\ref{irred}.

		Then there exists a number field $K$, a (possibly infinite) set $S$ of
		places of $K$, and a connected adjoint, absolutely simple $K$-group $\mathbf{H}$
		such that the following holds:

		Let $\tilde{\mathbf H}\to\mathbf H$ be the simply connected cover of $\mathbf H$ in the algebraic
		sense. Let $\calO(S)$ be the $S$-integers of $K$. The group $L\times D$
		embeds as a closed subgroup into the restricted (adelic) product
		$\prod'_{\nu\in S}\mathbf H(K_\nu)$.
		Under this embedding
		and under passing to a finite index subgroup, $\Gamma$ is contained
		in $\mathbf H(\calO(S))$ and
		the intersection of~$\Gamma$ with the image of $\prod'_{\nu\in S}\tilde{\mathbf H}(K_\nu)$
		is commensurable to the image of $\tilde{\mathbf H}(\calO(S))$.
\end{theorem}

The above theorem states essentially that all lattice inclusions $\Gamma<L\times D$ satisfying some natural conditions could be constructed from arithmetic data.
The following example describes, up to passing to finite index subgroups, all pairs of groups $\Gamma$, $D$ and embeddings $\Gamma<L \times D$
satisfying the condition of the theorem, for the special case $L=\PGL_2(\bbR)$.
These are obtained by the choice of the set $S$ and the subgroups $A$ and $B$ described below.
Similar classifications for other semi-simple groups $L$ could be achieved using 
Galois cohomology.

\begin{example}
Fix a possibly infinite set of primes $S\ne\emptyset$ and 
consider the localization $\bbZ_S<\bbQ$.
Fix a closed subgroup $A$ in the compact group $\prod_S \bbZ_p^\times/(\bbZ_p^\times)^2$
and a subgroup $B$ in the discrete group $\bigoplus_S \bbZ/2\bbZ$
(note that $\bbZ_p^\times/(\bbZ_p^\times)^2\simeq \bbZ/2\bbZ$ for $p>2$ and $\bbZ_2^\times/(\bbZ_2^\times)^2\simeq \bbZ/2\bbZ\times \bbZ/2\bbZ$).
The determinant homomorphism from $\prod'_S \GL_2(\bbQ_p)$ to the corresponding idele group,  which we naturally identify with $\prod_S \bbZ_p^\times \times \bigoplus_S \bbZ$, 
restricts to $\prod'_S \GL_2(\bbQ_p)\supset \GL_2(\bbZ_S)\to\bigoplus_S \bbZ<\prod_S \bbZ_p^\times \times \bigoplus_S \bbZ$. 
The determinant is well defined on $\PGL_2$ modulo squares, thus we obtain a map
$\prod'_S \PGL_2(\bbQ_p) \to \prod_S \bbZ_p^\times/(\bbZ_p^\times)^2 \times \bigoplus_S \bbZ/2\bbZ$
which restricts to $\PGL_2(\bbZ_S) \to \bigoplus_S \bbZ/2\bbZ$.
Let $D$ be the preimage of $A\times B$ under the first map and $\Gamma$ be the preimage of $B$ under the second.
Denote $L=\PGL_2(\bbR)$ and embed $\PGL_2(\bbZ_S)$ in $L$ via $\bbZ_S \hookrightarrow \bbR$. 
Thus $\Gamma$ embeds diagonally in $L\times D$
and its image is a lattice, whose projections to both $L$ and $D$ are dense and the projection of $\Gamma$ 
to $L$ is injective while $\Gamma$ satisfies~\ref{irred}.
\end{example}

\section{Sketch of proof of Theorem~\ref{thm: main}}\label{sec: ideas}

\subsubsection*{Step~1: Using~\ref{caf} to reduce to products.}

Burger and Monod~\cite[Theorem~3.3.3]{burger+monod}
observed that one obtains as a consequence of the positive
solution of Hilbert's 5th problem: Every locally compact group
has a finite index subgroup that
modulo its amenable radical splits as a product of a connected center-free
semi-simple real Lie group~$L$ without compact factors and a totally disconnected
group~$D$ with trivial amenable radical.
By Proposition~\ref{prop: compact amenable radical} we conclude that the amenable radical of $G$, thus that of
any of its finite index subgroups, is compact. 	
Therefore we may assume (up to passage to finite index subgroups and by dividing
out a normal compact subgroup) that $G$ is of the form $G=L\times D$ with $L$ and $D$ as
above.

If we assume that not all $\ell^2$-Betti numbers of $\Gamma$ vanish, we
can reach the same conclusion without appealing to
Proposition~\ref{prop: compact amenable radical} but by using
$\ell^2$-Betti numbers of locally compact groups~\cite{petersen} instead.
Since $\Gamma$ has a positive $\ell^2$-Betti number in some degree,
the same is true for $G$~\cite[Theorem~B]{kyed+petersen+vaes}, thus $G$ has a compact amenable radical~\cite[Theorem~C]{kyed+petersen+vaes}.

\subsubsection*{Step~2: Separating according to discrete and dense projections to the connected factor}

The connected Lie group factor $L$ splits as a product of
simple Lie groups $L=\prod_{i\in I} L_i$.
The projection of $\Gamma$ to $L$ might not have dense image. It is easy to see that
there is a maximal subset $J\subset I$ such that the projection $\pr_J$ of $\Gamma$ to
$L_J:=\prod_{j\in J}L_j$ has discrete image. Then $\Gamma_J=\pr_J(\Gamma)$
and $\Gamma'=\ker(\pr_J)\cap\Gamma$ are lattices in $L_J$ and $L_{J^c}\times D$,
respectively. So we obtain an extension of groups
\begin{equation}\label{eq: short exact sequence}
\Gamma'\hookrightarrow \Gamma\to \Gamma_J
\end{equation}
which are lattices in the
corresponding (split) extension of locally compact groups $L_{J^c}\times D\hookrightarrow L\times D\to L_J$. The projection of $\Gamma'$ to $L_{J^c}$ turns out to be dense.

Notice that finite generation of $\Gamma$ does not guarantee that $\Gamma_{J^c}$ is finitely generated. However $L_{J^c}\times D$ is still compactly
generated if $G$ is so.

\subsubsection*{Step~3: Distinguishing cases of the theorem}

Let $U<D$ be a compact open subgroup. Let $M:=\Gamma\cap (L\times U)$ and
$M':=\Gamma\cap (L_{J^c}\times U)$. We prove that $G$ is totally disconnected if
$M$ is finite and that $G=L$ if $M$ is infinite, but $M'$ is finite. The latter
step involves condition~\ref{irred}. Hence if $M$ or $M'$ is finite, the proof
is finished. In the remainder we discuss the case that $M'$ is infinite.
For simplicity let us first assume that the projection of $\Gamma'$ to $D$ is dense;
we return to this issue in the last step.

Consider $N':=\Gamma'\cap (\{1\}\times D)\triangleleft\Gamma'$ which can also be
regarded as a subgroup of $D$. As such it is also normal by denseness of the projection
$\Gamma'\to D$.
The assumptions of Theorem~\ref{thm: adelic arithmeticity} apart from condition~\ref{irred}
are satisfied for the lattice embedding $\Gamma'/N'<L_{J^c}\times D/N'$. From a more general
version of Theorem~\ref{thm: adelic arithmeticity} one concludes a posteriori that
$\Gamma'/N'$ satisfies~\ref{irred}, so the conclusion of Theorem~\ref{thm: adelic arithmeticity}
holds true for $\Gamma'/N'<L_{J^c}\times D/N'$. Because of
compact generation of $L_{J^c}\times D/N'$ we can exclude the adelic case
and conclude that $\Gamma'/N'$ is an $S$-arithmetic lattice for a finite set $S$ of
primes.

\subsubsection*{Step~4: Using~\ref{nbc} to identify group extensions}

In this step we show that $N'$ is finite, thus trivial (since $D$ has no compact
normal subgroups). The proof involves the use of~\ref{nbc} for $N'$ and $M'$ and
Margulis' normal subgroup theorem. Hence $\Gamma'$ is an $S$-arithmetic lattice. As
such it has a finite outer automorphism group which implies that, after passing
to finite index subgroups, the extension~(\ref{eq: short exact sequence}) splits.
By~\ref{irred} $\Gamma<G$ is an $S$-arithmetic lattice embedding.

\subsubsection*{Step~5: Quasi-isometric rigidity results}

We have previously assumed that the projection $\Gamma'\to D$ is dense. If it is not
we have to identify the difference between the closure $D'$ of the image of the projection
and $D$. The subgroup $D'<D$ is cocompact, thus $D'\hookrightarrow D$
is a quasi-isometry. By the argument before we know that $D'$ is a product of algebraic
groups over non-Archimedean fields and thus acts by isometries on a product $B$ of
Bruhat-Tits buildings. 
By conjugating with $D'\hookrightarrow D$ we obtain a homomorphism
of $D'$ to the quasi-isometry group of $B$. We finally appeal to the quasi-isometric rigidity
results of Kleiner-Leeb~\cite{kleiner+leeb} and Mosher-Sageev-White~\cite{mosher+sageev+whyte} to conclude
that $\Gamma<G$ is an $S$-arithmetic lattice up to tree extension.
\section{Special cases}

It is instructive to investigate the consequences of Theorem~\ref{thm: main}
for specific groups. Rather than
just applying Theorem~\ref{thm: main} we sketch a blend of ad hoc arguments and techniques
of the proof of Theorem~\ref{thm: main} to most easily classify all lattice embeddings in the
following three cases.

\subsection{$\Gamma$ is a free group.}
Let $\Gamma$ be a non-commutative finitely generated free group. Let $\Gamma<G$ be a lattice
embedding. We show that, up to finite index and dividing out a normal compact subgroup,
$F<G$ is $\PSL_2(\bbZ)<\PSL_2(\bbR)$ or $G$ embeds as a closed cocompact subgroup
in the automorphism group of a tree.

As explained in the first step of Subsection~\ref{sec: ideas} one
can avoid the use of Proposition~\ref{prop: compact amenable radical}
by using the positivity of the first $\ell^2$-Betti number of $\Gamma$
to conclude that $G$ has a compact amenable radical. Up to passage to a finite index subgroup
and dividing out a compact amenable radical we may assume that $G$ is a product
$G\cong L\times D$. By the
K\"unneth formula $L\times D$ can have positive first $\ell^2$-Betti number only if
one of the factors is compact. Since $G$ has trivial amenable radical,
this implies that $G$ is either $L$ or $D$. In the first case $G$ must be $\PSL_2(\bbR)$.
In the second case $G$ is totally disconnected, and $\Gamma<G$ as a torsion-free lattice must
be cocompact. By~\cite[Theorem~9]{mosher+sageev+whyte} $G$ embeds as a closed cocompact
subgroup of the automorphism group of a tree.

\subsection{$\Gamma$ is a surface group.}

Let $\Gamma$ be the fundamental group of a closed oriented surface~$\Sigma_g$ of genus~$g\ge 2$.
Let $\Gamma<G$ be a lattice embedding.
Similarly as for free groups, by using the positivity of the first $\ell^2$-Betti number, we
conclude that $G$, up to passage to a finite index subgroup
and dividing out a compact amenable radical, is either $\PSL_2(\bbR)$ or a totally disconnected
group with trivial amenable radical. In the latter case $\Gamma$ is cocompact.

We argue that the totally disconnected case cannot happen unless $G$ is discrete and so $\Gamma<G$
is the trivial lattice embedding: The inclusion $\Gamma\to G$
is a quasi-isometry in that case. So we obtain a
homomorphism $G\to \QI(G)\cong \QI(\Gamma)$. Each quasi-isometry
induces a homeomorphism of the boundary $\partial\Gamma\cong S^1$, so we obtain a homomorphism
$f: G\to \QI(\Gamma)\to \Homeo_+(S^1)$. One can verify that $f$ is continuous~\cite[Theorem~3.5]{furman-mostow}
and $\ker(f)$ is compact, thus trivial
by the triviality of the amenable radical. Let $U<G$ be a compact-open subgroup. Then $f(U)<\Homeo_+(S^1)$ is
a compact subgroup, hence $f(U)$ is either finite or isomorphic to $\SO(2)$~\cite[Lemma~3.6]{furman-mostow}.
But it cannot be isomorphic to a connected group. Therefore $f(U)$ is finite, which implies that $G$ is discrete.

%

\subsection{$\Gamma=\PSL_n(\bbZ[1/p])$, $n\ge 3$.} 

Recall that $\Gamma$ embeds as a non-uniform lattice in $\PSL_n(\bbR)\times\PSL_n(\bbQ_p)$ via $\bbZ[1/p]\to\bbR\times \bbQ_p$;
we denote by $\pr_1:\Gamma\to\PSL_n(\bbR)$ and $\pr_2:\Gamma\to\PSL_n(\bbQ_p)$ the injective projections.
Let us verify~\ref{caf}:
For any commensurated amenable subgroup $A<\Gamma$, the connected component $H^0$ of the Zariski closure 
$H=\overline{\pr_1(A)}$ is amenable and normal in $\PSL_n(\bbR)$ because replacing $A$ 
by a finite index subgroup does not change $H^0$. 
Hence $H^0$ is trivial, and so $H$ and $A$ are finite. 

Let $\Gamma$ be embedded as a lattice in some locally compact group $G$.
Using \ref{caf} as in the first step of Subsection~\ref{sec: ideas} we replace $G$ by $L\times D$,
where $L$ is a (possibly trivial) connected real Lie group, 
$D$ is totally disconnected, and both have trivial amenable radicals.
Let $E<D$ denote the closure of $\pr_D(\Gamma)$; then $\Gamma<L\times E$ where $E$ is totally disconnected,
and $E$ has finite covolume in $D$.

%
%
Case $L=\{1\}$ corresponds to the trivial lattice embedding $\Gamma<\Gamma$. 
Indeed, in this case $\Gamma$ is a lattice in a totally disconnected $D$, and having bounded torsion, it is cocompact.
This allows us to use the results on quasi-isometric rigidity \cite{wortman} and obtain a homomorphism 
$D\to \QI(D)= \QI(\Gamma)\simeq \PSL_n(\bbQ)$ (hereafter $\simeq$ stands for commensurability) 
that can be further shown to have an image commensurable to $\Gamma$.

If $L$ is non-trivial, then it is a center-free, semi-simple Lie group without compact factors.
By Borel's density theorem the projection map $\pr_L:\Gamma\to L$ has Zariski dense image, 
and Margulis' superrigidity implies that $L= \PSL_n(\bbR)$ and $\pr_L=\pr_1$ 
(this can also be shown by elementary means by conjugating unipotent matrices). 
Let $E<D$ denote the closure of $\pr_D(\Gamma)$; we get a lattice embedding $\Gamma<L\times E$ where $E$ is totally disconnected,
and $E$ has finite covolume in $D$.
If $U<E$ is a compact open subgroup, $\Gamma_U=\Gamma\cap (\PSL_n(\bbR)\times U)$ is a lattice in $\PSL_n(\bbR)\times U$,
that projects to a lattice $\Delta<\PSL_n(\bbR)$.
We claim that $\Delta\simeq \PSL_n(\bbZ)$.
Indeed, it follows from Margulis' superrigidity (recall that $n\ge 3$) 
applied to $\pr_2\circ \pr_1^{-1}:\Delta\to \Gamma_U\to \PSL_n(\bbQ_p)$ 
that $\pr_2(\Gamma_U)$ is contained in a maximal compact subgroup commensurated to $\PSL_n(\bbZ_p)$, yielding $\Gamma_U\simeq \PSL_n(\bbZ)$.

Let $H$ be the closure in $E\times \PSL_n(\bbQ_p)$ of 
$\Lambda=\{(\pr_E(\gamma),\pr_2(\gamma)) \mid \gamma\in\Gamma\}$. Then $H\cap (U\times \PSL_n(\bbQ_p))$ 
is the closure of $\Lambda\cap (U\times \PSL_n(\bbQ_p))$, thus compact because of $\Gamma_U\simeq \PSL_n(\bbZ)$. 
Also $\Ker(\pr_{E}:H\to E)$ is compact. The group $\pr_E(H\cap (U\times \PSL_n(\bbQ_p)))$ is compact, hence closed 
and equals $\overline{\pr_E(\Lambda\cap (U\times\PSL_n(\bbQ_p))}=U$. Since $\pr_E(H)$ is dense in $E$ and 
contains the open subgroup~$U$, we obtain that $\pr_E(H)=E$. Similarly, $\pr_2(H)=\PSL_n(\bbQ_p)$.
This implies that $H$ is a graph of a continuous surjective homomorphism $E\to \PSL_n(\bbQ_p)$ whose 
kernel $K$ is contained in $U$, thus compact.

Finally, to recover the original group $D$, one uses the fact that $E$ is cocompact in $D$, so 
quasi-isometric rigidity~\cite{kleiner+leeb} of the Bruhat-Tits building $X_n$ of $\PSL_n(\bbQ_p)\cong E/K$ gives 
$D\to \QI(D)\cong \QI(E)\cong \QI(X_n)\simeq \PSL_n(\bbQ_p)$.

\bibliographystyle{abbrv}
\bibliography{latenv-short-arx}

\end{document}